\title{On Vervaat transform of Brownian bridges and Brownian motion}
\author{Jim Pitman \thanks{Statistics department, University of California, Berkeley.  Email: pitman@stat.berkeley.edu} and Wenpin Tang \thanks{ Statistics department, University of California, Berkeley. Email: wenpintang@stat.berkeley.edu.} \footnote{The author is also affiliated to Ecole Polytechnique during the preparation of the work.}}

\date{\today}
\documentclass[12pt,a4paper]{article}
\usepackage{graphicx} \usepackage{amsmath, amsfonts}
\usepackage{float}
\usepackage[utf8]{inputenc} \usepackage[T1]{fontenc} \usepackage[french,english]{babel} \usepackage{color}
\usepackage[top=1in, bottom=1in, left=1in, right=1in]{geometry}
\DeclareMathOperator\erf{erf}
\DeclareMathOperator\argmin{argmin}
\newtheorem{theorem}{Theorem}[section]
\newtheorem{lemma}[theorem]{Lemma}
\newtheorem{proposition}[theorem]{Proposition}
\newtheorem{corollary}[theorem]{Corollary}
\newtheorem{conjecture}[theorem]{Conjecture}
\begin{document}
\maketitle
\textbf{Abstract:} For a continuous function $f \in \mathcal{C}([0,1])$, define the Vervaat transform $V(f)(t):=f(\tau(f)+t \mod1)+f(1)1_{\{t+\tau(f) \geq 1\}}-f(\tau(f))$, where $\tau(f)$ corresponds to the first time at which the minimum of $f$ is attained. Motivated by recent study of quantile transforms of random walks and Brownian motion, we study the Vervaat transform of Brownian motion and Brownian bridges with arbitary endpoints. When the two endpoints of the bridge are not the same, the Vervaat transform is not Markovian. We describe its distribution by path decomposition and study its semimartingale property. The expectation and variance of the Vervaat transform of Brownian motion are also derived.



\setcounter{tocdepth}{1}
\tableofcontents
\section{Introduction}
In a recent work of Assaf et al \cite{AFP}, a novel path transform, called the quantile transform $Q$ has been studied both in discrete and continuous time settings. 
Inspired by previous work in fluctuation theory (see e.g. and Wendel \cite{Wendel} and Port \cite{Port}), the quantile transform for simple random walks is defined as follows. For $w$ a simple walk of length $n$, with increments of $\pm 1$
the quantile transform associated to $w$ is defined by:
$$\forall j \in [1,n],~Q(w)_j:=\sum_{i=1}^j w(\phi_w(i))-w(\phi_w(i)-1),$$
where $\phi_w$ is the quantile permutation on $[1,n]$ 
defined by lexicographic ordering on pairs $(w(j-1),j)$, that is $w(\phi_w(i)-1)<w(\phi_w(j)-1)$ or $w(\phi_w(i)-1)=w(\phi_w(j)-1)$, $\phi_w(i) \leq \phi_w(j)$ if and only if $i \leq j$.

As shown in \cite{AFP}, the scaling limit of this transformation of simple random walks is the quantile transform in the continuous case of Brownian motion $B:= (B_t; 0 \leq t \leq 1)$:
$$\forall t \in [0,1],~Q(B)_t:=\frac{1}{2} L_1^{a(t)} + (a(t))^{+}-(a(t)-B_1)^{+},$$
where $L_1^a$ is the local time of $B$ at level $a$ up to time $1$ and $a(t):=\inf\{a; \int_0^1 1_{B_{s} \leq a}ds>t\}$ is the quantile function of occupation measure (see Dassios \cite{Dassios}, Embrechts et al \cite{ERY} for general background).

The key result of  \cite{AFP} is to identify the distribution of the somewhat mysterious $Q(w)$ with that of the Vervaat transform $V(w)$ defined as:
$$V(w)_i: =  \left\{ \begin{array}{ccl} w_{\tau_n+i}-w_{\tau_n}  &\mbox{for} & i \leq n-\tau_n \\ w_{\tau_n+i-n}+w_n-w_{\tau_n} &\mbox{for} & n-\tau_n \leq i \leq  n,  \end{array}\right.$$
where $\tau_n:= \argmin_{i \in [1,n]} w_i$. Consequently, $(Q(B); 0 \leq t \leq 1) \stackrel{d}{=} (V(B); 0\leq t \leq 1)$ for 
$$V(B)_t:=  \left\{ \begin{array}{ccl} B(1-A+t)-B(1-A)  &\mbox{for} & 0 \leq t \leq A \\ B(t-A)+B(1)-B(1-A) & \mbox{for} & A \leq t \leq 1,  \end{array}\right.$$
where $A$ is the a.s. arcsine split ($1-A:=\argmin_{t \in [0,1]} B_t$, see Karatzas and Shreve \cite{KS}).

As a result, to understand quantile transform of $B$, it is equivalent to study its substitute, the Vervaat transform $V(B)$. Historically, Vervaat \cite{Vervaat} showed that if $B$ is conditioned to both start and
end at $0$, then $V(B)$ is a Brownian excursion:
\begin{theorem}\cite{Vervaat}
\label{Vervaat}
$(V(B^{0,br}); 0 \leq t \leq 1) \stackrel{d}{=} (B^{ex}; 0 \leq t \leq 1)$, where $(B^{0,br}_t; 0 \leq t \leq 1)$ is a Brownian bridge of length $1$ starting at $0$ and ending at $0$.
\end{theorem}
Biane \cite{Biane} proved a converse theorem to Vervaat's result, i.e. recover standard Brownian bridges from Brownian excursion by uniform sampling:
\begin{theorem}\cite{Biane}
\label{Biane}
 Let $B^{ex}$ be a standard Brownian excursion and $U$ a uniformly distributed random variable independent of $B^{ex}$. Then the shifted process $\theta(B^{ex},U)$ defined by
$$\theta(B^{ex},U)_t: =  \left\{ \begin{array}{ccl} B^{ex}_{U+t}-B^{ex}_U  &\mbox{for} & 0 \leq t \leq 1-U \\ B^{ex}_{U+t-1}-B^{ex}_U &\mbox{for} & 1-U \leq t \leq 1,  \end{array}\right.$$
is a standard Brownian bridge.
\end{theorem}
Chaumont \cite{Chaumont} extended partly the result to stable cases, Chassaing and Jason \cite{CJ} to the reflected Brownian bridges case, Miermont \cite{Miermont} to the spectrally positive case, Fourati \cite{Fourati} to the general L\'{e}vy case under some mild hypotheses, Le Gall and Weill \cite{LeGall} to the Brownian tree case and more recently, Lupu \cite{Lupu} to the diffusion case. However, as far as we are aware, there has not been previous study of the Vervaat transform of an unconditioned Brownian motion $B$ or of the Brownian bridges $B^{\lambda,br}$ ending at $\lambda \ne 0$.

The contribution of the current paper is to give some path decomposition result of Vervaat transform of Brownian bridges (for simplicity, call them Vervaat bridges) with non-zero endpoints. In the case of a Vervaat bridge with negative endpoint $V(B^{\lambda,br})$ where $\lambda<0$, the key idea is to decompose it into two pieces, the first piece a Brownian excursion and the second piece a first passage bridge. The main result is stated as follows:
 \begin{theorem} 
\label{PT}
Let $\lambda<0$. Given $Z^{\lambda}$ the first return to $0$ of $V(B^{\lambda,br})$, whose density is given by
\begin{equation}
\label{aldous}
f_{Z^{\lambda}}(t)=\frac{|\lambda|}{\sqrt{2 \pi t(1-t)^3}} \exp\left(-\frac{\lambda^2 t}{2(1-t)}\right),
\end{equation}
the path is decomposed into two (conditionally) independent pieces:
\begin{itemize}
\item
$(V(B^{\lambda,br})_{u}; 0 \leq u \leq Z^{\lambda})$ is a Brownian excursion of length $Z^{\lambda}$;
\item
 $(V(B^{\lambda,br})_{u}; Z^{\lambda} \leq u \leq 1)$ is a first passage bridge through level $\lambda$ of length $1-Z^{\lambda}$.
\end{itemize}
\end{theorem}
\begin{center}
\includegraphics[width=0.6\textwidth]{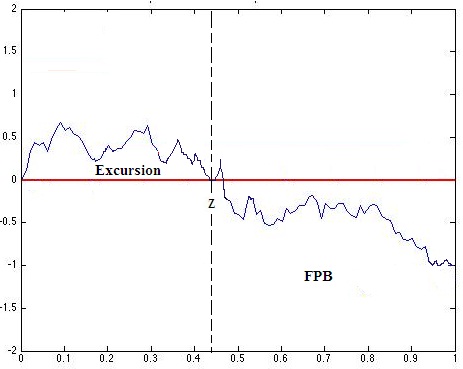}\\
Fig 1. Vervaat bridge $=$ Excursion $+$ First passage bridge.
\end{center}
Note that Theorem \ref{Vervaat} \cite{Vervaat} is recovered as a weak limit $\lambda \rightarrow 0$ of the previous theorem.
The parametric density family $(f_{Z^{\lambda}})_{\lambda<0}$ appears earlier in the work of Aldous and Pitman \cite{AP}, Corollary $5$ when they studied the standard additive coalescent. Precisely, $Z^{\lambda} \stackrel{d}{=}\frac{B_1^2}{\lambda^2+B_1^2}$ where $B_1$ is normal distributed with mean $0$ and variance $1$. We also refer readers to Pitman \cite{Pitman}, Chapter 4 for some discussion therein.

For the Vervaat bridges $V(B^{\lambda,br})$ which ends up with some positive value, it is easy to see that we have the following duality relation:
\begin{equation}
\left(V(B^{\lambda,br})_t; 0 \leq t \leq 1\right) \stackrel{d}{=} \left(V(B^{-\lambda,br})_{1-t}+\lambda; 0 \leq t \leq 1\right) \quad \mbox{for}~\lambda>0.
\end{equation}
In other words, looking backwards, we have a first piece of excursion above level $\lambda$ followed by a first passage bridge. Note in addition that a first passage bridge form $\lambda>0$ to $0$ has the same distribution as a three dimensional Bessel bridge from $\lambda$ to $0$ (see Biane and Yor \cite{BY}). We have the following decomposition of Vervaat bridges with negative endpoints:
\begin{corollary}
\label{PT2}
Let $\lambda>0$. Given $\widehat{Z}^{\lambda}$ the time of last hit of $\lambda$ by $V(B^{\lambda,br})$ strictly before $1$, whose density is given by $f_{\widehat{Z}^{\lambda}}(t)=f_{Z^{-\lambda}}(1-t)$ as in \eqref{aldous}, the path is decomposed into two (conditionally) independent pieces:
\begin{itemize}
\item
 $(V(B^{\lambda,br})_{u};0 \leq u \leq \widehat{Z}^{\lambda})$ is a three dimensional Bessel bridge of length $\widehat{Z}^{\lambda}$ starting from $0$ and ending at $\lambda$;
\item
$(V(B^{\lambda,br})_{u};\widehat{Z}^{\lambda} \leq u \leq 1)$ is a Brownian excursion above level $\lambda$ of length $1-\widehat{Z}^{\lambda}$. 
\end{itemize}
\end{corollary}
\begin{center}
\includegraphics[width=0.6\textwidth]{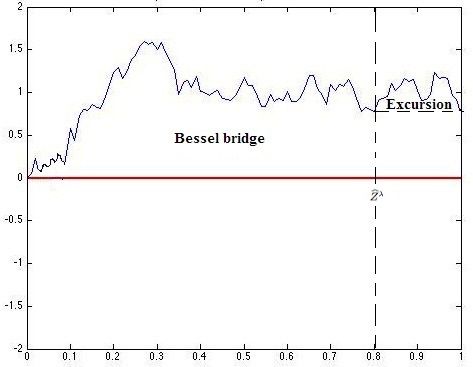}\\
Fig 2. Vervaat bridge$=$Bessel bridge$+$Excursion.
\end{center}

The rest of the paper is organized as follows. In Section 2, we provide the proof of Theorem \ref{PT} by random walk approximation, which is based on some bijection lemma proved in Assaf et al\cite{AFP}.

In Section 3, we give a thorough study of $V(B^{\lambda,br})$ where $\lambda \neq 0$ using Theorem \ref{PT} and Corollary \ref{PT2}. We prove that such processes are not Markov (Section 3.2). We also relate these processes to some simpler ones (Section 3.1, 3.3) and study the convex minorant of such processes (Section 3.4). 

In Section 4, we focus on studying the Vervaat transform of Brownian motion. We first prove that $V(B)$ is not Markov as well (Section 4.1). Nevertheless, we show that it is a semimartingale using Bichteler-Dellacherie's characterization for semimartingales (Section 4.2). Finally, we provide explicit formulae for the first two moments of the Vervaat transform of Brownian motion (Section 4.3).
\section{Path decomposition for Vervaat bridges}
The whole section is devoted to proving Theorem \ref{PT}. We use a discrete approximation argument to obtain the path decomposition of $V(B^{\lambda,br})$ where $\lambda<0$. Also we obtain an analog to Theorem \ref{Biane} as a by-product.
\subsection{Discrete case analysis} 
We begin with the discrete time analysis of random walk cases which is based on combinatorial principles. For a simple random walk $w$ of length $n$ with increments $\pm 1$, we would like to describe the law of $V(w^{a}):=(V(w)|w(n)=a)$ where $a<0$ having the same parity as $n$.

Denote $\tau_V(w)=\min\{j \in [0,n]; w(j) \leq w(i), \forall i \in [0,n] \}$ (the first global minimum of the path) and $K(w)=n-\tau_V(w)$ (distance from the first global minimum to the end of the path). Following from Theorem $7.3$ in Assaf et al \cite{AFP}, the mapping $w \rightarrow (V(w), K(w))$ is a bijection between $walk(n)$, the set of simple random walks of length $n$ and the set 
$$\{(v,k); v \in walk(n), v(j) \geq 0 ~\mbox{for} ~0 \leq j\leq k ~\mbox{and}~v(j)>v(n)~\mbox{for}~ k \leq j <n\},$$
where $k$, called a {\em helper variable}, records the splitting position in the original path.

The following result turns out to be a direct consequence of this theorem related to Vervaat bridges.
\begin{lemma} 
\label{FP}
$w^a \rightarrow (V(w^a), K(w^a))$ forms a bijection between $\{w \in walk(n): w(n)=a\}$ (simple random walk bridges which end at $a<0$) and the set 
$$\{(v,k); v \in walk(n), v(j) \geq 0 ~\mbox{for} ~0 \leq j\leq k ,~v(j)>a~\mbox{for}~ k \leq j <n ~\mbox{and}~v(n)=a\}.$$
\end{lemma}
Observe that, to each pair $(v,k)$ in the above set, one can associate a unique triple $(Z^{a}, f_{Z^{a}}^{br,1},f_{Z^{a}}^{br,2})$ where 
\begin{itemize}
\item $Z^{a}$ is the first time that the path hits level $-1$,
\item  $f_{Z^{a}}^{br,1}$ is the sample path of a first passage bridge of length $Z^{a}$ through level $-1$,
\item $f_{Z^{a}}^{br,2}$ is that of a first passage bridge of length $n-Z^{a}$ starting at $-1$ through $a$. 
\end{itemize}
Remark that to different pairs $(v,k)$, one may have the same triple $(Z^{a}, f_{Z^{a}}^{br,1},f_{Z^{a}}^{br,2})$.

We now focus on calculating explicitly the distribution of $Z^{a}$ by counting paths. By Lemma \ref{FP}, the total number of the Vervaat transform paths (counting with multiplicity) is $\binom{\frac{n+|a|}{2}}{n}$  since the mapping $w \rightarrow V(w)$ is not injective.

Moreover, the number of first passage bridges through level $-1$ of odd length $l$ is $\frac{1}{l} \binom{\frac{l+1}{2}}{l}$ and the number of first passage bridges starting at $-1$ through level $a$ of length $n-l$ is $\frac{|a|-1}{n-l} \binom{\frac{n-l+|a|-1}{2}}{n-l}$ (see Chapter III of Feller \cite{Feller}). Therefore, the total number of the Vervaat transform configurations (counting with multiplicity ) is 
$$\frac{|a|-1}{l(n-l)} \binom{\frac{l+1}{2}}{l} \binom{\frac{n-l+|a|-1}{2}}{n-l}.$$
Also note that every Vervaat transform configuration is counted exactly $l$ times (by bijection lemma \ref{FP}). Hence,
\begin{equation}
\label{1}
\mathbb{P}(Z^{a}=l)=\frac{|a|-1}{n-l} \frac{\binom{\frac{l+1}{2}}{l} \binom{\frac{n-l+|a|-1}{2}}{n-l}}{\binom{\frac{n+|a|}{2}}{n}}.
\end{equation}

Combining the above discussions, we get the following path decomposition result for discrete Vervaat bridges with negative endpoint:
\begin{theorem} 
\label{discrete}
Let $a<0$ and have the same parity as $n$. Given $Z^{a}:=\min\{j>0;V(w^{a})_j=-1\}$ (distributed as \eqref{1}), the path is decomposed into two (conditionally) independent pieces:
\begin{itemize}
\item
$V(w^{a})|_{[0,Z^{a}]}$ is a random walk first passage bridge of length $Z^{a}$ through level $-1$,
\item
$V(w^{a})|_{[Z^{a},n]}$ is a random walk first passage bridge starting at $-1$ through level $a$ of length $n-Z^{a}$.
\end{itemize}
\end{theorem}
\begin{center}
\includegraphics[width=0.6\textwidth]{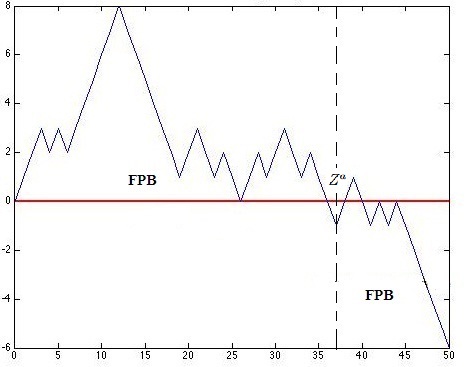}\\
Fig 3. Discrete Vervaat bridge$=$First passage bridge$+$First passage bridge.
\end{center}

The theorem provides a path decomposition of Vervaat bridges into two pieces of first passage bridges, one through level $-1$ and the other from $-1$ to $a$. Note that it is also possible to decompose the path slightly differently by a first piece of excursion and the second a first passage bridge through level $-a$. However, the distribution of the splitting position is much less explicit and thus does not make the proof any easier when passing to the scaling limit.
\subsection{Continuous case: passage to weak limit}
We now turn to the continuous case by appealing to invariance principles. We derive the path decomposition result from Theorem \ref{discrete}. 

For $\lambda<0$ and $0<t<1$, let $\lambda_n \sim \lambda \sqrt{n}$ and have the same parity as $n$ and $t_n := 2[\frac{tn}{2}]+1$ be two fixed sequences. Let $S^{\lambda_n}$ be simple random walks of length $n$ with increments $\pm 1$ which end at $\lambda_n$, $V(S^{\lambda_n})$ be the associated discrete Vervaat bridge and $Z^{\lambda_n}:= \inf \{j>0: V(S^{\lambda}_n)_j=-1\}$. Define $\left(V(S^{\lambda_n})(u); 0 \leq u \leq n \right)$ to be the linear interpolation of the discrete Vervaat bridge $V(S^{\lambda_n})$. 

Recall some invariance principle results on metric space $\mathcal{C}[0,1]$ (continuous functions on $[0,1]$). For general background on weak convergence in $\mathcal{C}[0,1]$, we refer the readers to Chapter 2, Billingsley \cite{Bill2}. 
\begin{lemma}
(a).$(\frac{1}{\sqrt{n}}V(S^{\lambda_n})(nu); 0 \leq u \leq 1)$ converges in $\mathcal{C}[0,1]$ to $(V(B^{\lambda,br})_u; 0 \leq u \leq 1)$. (b).Given $Z^{\lambda_n}=t_n$, $(\frac{1}{\sqrt{n}}V(S^{\lambda_n})(nu); 0 \leq u \leq t)$ converges in $\mathcal{C}[0,1]$ to a Brownian excursion of length $t$ and $(\frac{1}{\sqrt{n}}V(S^{\lambda_n})(nu); t \leq u \leq 1)$ converges in $\mathcal{C}[0,1]$ to a first passage bridge through level $ \lambda$ of length $1-t$, (conditionally) independent of the excursion.
\end{lemma}
\textbf{Proof:} The assertion $(a)$ can be viewed as a variant of the results proved in Vervaat \cite{Vervaat}. According to Theorem \ref{discrete}, given $Z^{\lambda_n}=t_n$, the path of $V(S^{\lambda_n})$ is split into two (conditionally) independent pieces of discrete first passage bridges. Following  Bertoin et al \cite{BCP} and Iglehart \cite{Iglehart}, the scaled first passage bridge through level $-1$ converges weakly to a Brownian excursion and the scaled first passage bridge from $-1$ to $\lambda_n$ converges weakly to a first passage bridge through level $\lambda$. This proves $(b)$.  $\square$

To prove Theorem \ref{PT}, we need to compute the limiting distribution of $Z^{\lambda_n}=t_n$ as $n \rightarrow \infty$. Precisely,
\begin{equation}
\label{3}
n \mathbb{P}(Z^{{\lambda}_n}=t_n)= \frac{n|\lambda_n|}{n-t_n} \frac{\binom{\frac{t_n+1}{2}}{t_n} \binom{\frac{n-t_n+|\lambda_n|-1}{2}}{n-t_n} }{\binom{\frac{n+|\lambda_n|}{2}}{n}}.
\end{equation}
Using Stirling's formula, we see:
$$\binom{\frac{t_n+1}{2}}{t_n} \sim \sqrt{\frac{2}{\pi nt}} 2^{nt};$$
$$ \binom{\frac{n+|\lambda_n|}{2}}{n} \sim \sqrt{\frac{2}{\pi n}} 2^n \exp \left(-\frac{\lambda^2}{2}\right);$$
and
$$\binom{\frac{n-t_n+|\lambda_n|-1}{2}}{n-t_n} \sim \sqrt{\frac{2}{\pi n(1-t)}} 2^{n(1-t)} \exp\left(-\frac{\lambda^2}{2(1-t)}\right).$$
Injecting these terms in \eqref{3}, we deduce the limiting distribution as $n \rightarrow \infty$ given by \eqref{aldous} . By a  local limit argument (see Billingsley \cite{Bill}, Exercise $25.10$), we conclude that $Z^{\lambda}$ has density $f_{Z^{\lambda}}$ given in \eqref{aldous}.\\\\
\textbf{Remark:} P.Fitzsimmons points out that the decomposition result is also a consequence of a local Williams decomposition, which can be found in the section $6$ of \cite{Fitz}.

The next theorem is a direct consequence of Theorem \ref{PT} and should be called a corollary at best. Because of its importance, however, we give it status of a theorem.
\begin{theorem}
Given $Z^{\lambda}$ the length of first excursion of $(V(B^{\lambda,br})_t; 0 \leq t \leq 1)$ where $\lambda<0$, the split position $A^{\lambda}:=1-\argmin_{t \in [0,1]} B^{\lambda,br}_t$ (distance from the minimum of the original bridge path to the end) is (conditionally) independent of $V(B^{\lambda,br})$ and uniformly distributed on $[0,Z^{\lambda}]$, In particular, its density is
$$f_{A^{\lambda}}(a)=\int_a^1 \frac{f_{Z^{\lambda}}(t)}{t}dt,$$
where $f_{Z^{\lambda}}$ is given by \eqref{aldous}.
\end{theorem}
\textbf{Proof:} Note that in the discrete case, given a Vervaat bridge path, the {\em helper variable} $k$ takes values exactly in $\{0,..., Z^{\lambda_n}\}$ where $Z^{\lambda_n}$ is the first time that the path returns to $0$. This implies that given $Z^{\lambda_n}$, the minimum position of the original bridge is uniformly distributed on $[0,Z^{\lambda_n}]$. We then obtain the results in the theorem by passing to the scaling limit.  $\square$
\begin{corollary}
Let $V(B^{\lambda,br})$ be the Vervaat transform of a Brownian bridge ending at $\lambda<0$. Given $Z^{\lambda}$ the first return to $0$ of $V(B^{\lambda,br})$, let $A^{\lambda}$ be uniformly distributed on $[0,Z^{\lambda}]$. Then the shifted process $\theta(V(B^{\lambda,br}),A^{\lambda})$ as defined in Theorem \ref{Biane} is a Brownian bridge ending at $\lambda$ which attains its minimum at $1-A^{\lambda}$.
\end{corollary}
\textbf{Remark:} The above corollary holds true for $\lambda \leq 0$ and the case $\lambda=0$, i.e. Theorem \ref{Biane} \cite{Biane} is recovered as a weak limit $\lambda \rightarrow 0$: $Z^{\lambda} \stackrel{d}{\rightarrow} 1$ and $A^{\lambda} \stackrel{d}{\rightarrow} \mbox{Uniform}[0,1]$.
\section{Study of Vervaat bridges}
In this section, we will study thoroughly the Vervaat bridges with non-zero endpoint. First, we give an alternative construction of $V(B^{\lambda,br})$ using length-biased sampling techniques. Next we show that such processes are not Markov with respect to their induced filtrations. Despite lack of markovianity, they are semimartingales with respect to their own filtrations and the proof of this fact is reported to Section $4.2$. Moreover, we relate Vervaat bridges to drifting excursion by additional conditioning. To close the section, we study some properties of convex minorant of $V(B^{\lambda,br})$ where $\lambda<0$.
\subsection{Construction of Vervaat bridges via Brownian bridges}
In the current part, we try to provide an alternative construction of the Vervaat bridges with negative endpoint via standard Brownian bridges (which end at $0$). It is obvious that the Vervaat bridges with positive endpoint can be treated similarly by time reversal.

Let $\lambda<0$. As seen in the last section, conditioned on $Z^{\lambda}$ the first return to $0$, the process is split into $B^{ex,Z^{\lambda}}$ an excursion of length $Z^{\lambda}$ followed by $F^{\lambda,1-Z^{\lambda}}$ a first passage bridge through $\lambda$ of length $1-Z^{\lambda}$, independent of each other. Formally, $V(B^{\lambda,br})$ looks much like a standard first passage bridge (of length 1) except that it has an excursion piece placed first. Therefore, it is interesting to ask whether this process can be derived from standard first passage bridge via some simple operations.

Recall that a standard first passage bridge can be constructed via standard Brownian bridge by conditioning on its local time. Denote $(F^{\lambda}_t; 0 \leq t \leq 1)$ for a standard first passage bridge through $\lambda<0$. Following from Bertoin et al \cite{BCP}, 
\begin{equation}
\label{4}
 (F^{\lambda}_t; 0 \leq t \leq 1) \stackrel{d}{=} (|B^{0,br}_t|-L_t^0(B^{0,br}); 0 \leq t \leq 1|L_1^0(B^{0,br})=|\lambda|),
\end{equation}
where $L_t^0$ is the local time (of a Brownian bridge) at level $0$ up to time $t$. 

In light of the above construction, the following theorem tells how to construct the Vervaat bridges with negative terminal value by standard Brownian bridges.
\begin{theorem}
Let $U$ be uniformly distributed on $(0,1)$ independent of  $X:=(B^{0,br}(t); 0 \leq t \leq 1|L_1^0(B^{0,br})=|\lambda|)$ where $\lambda<0$ and $(G_U,D_U)$ be the signed excursion interval which contains $U$. Let $\tilde{X}$ be the process by exchanging the position of the excursion of $X$ straddling time $U$ and the path along $[0,G_U]$, namely:
$$\tilde{X}_t =     \left\{ \begin{array}{rcl}  X_{t+G_U} & \mbox{for} & 0 \leq t \leq D_U-G_U \\ X_{t-D_U+G_U}  & \mbox{for} & D_U-G_U \leq t \leq D_U \\  X_t & \mbox{for} &  D_U \leq t \leq 1.\end{array}\right.$$
Then we have the following identity in law:
\begin{equation}
\label{5}
(|\tilde{X}(t)|-L_t^0(\tilde{X}); 0 \leq t \leq 1|L_1^0(\tilde{X})=|\lambda|) \stackrel{d}{=} (V(B^{\lambda,br})_t; 0 \leq t \leq 1).
\end{equation}
\end{theorem}
\begin{center}
\includegraphics[width=0.9\textwidth]{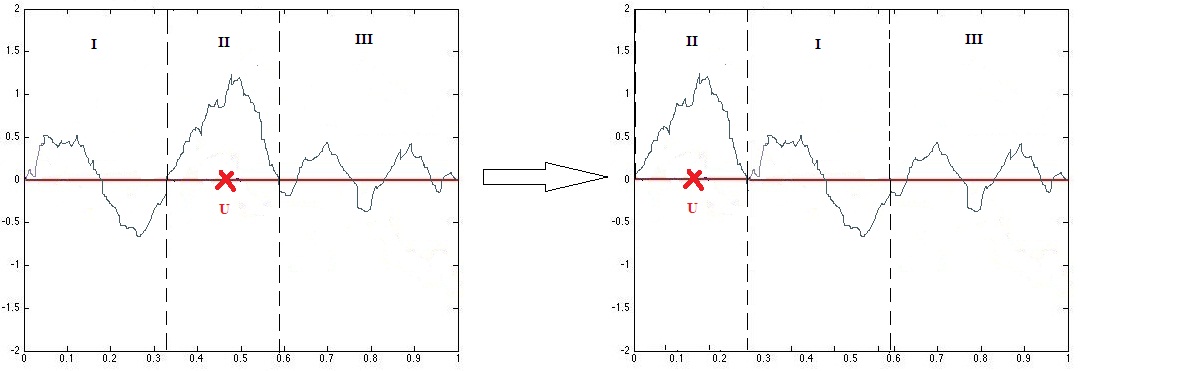}\\
Fig 4. Length-biased pick for a Brownian bridge conditioned on its local time.
\end{center}
\textbf{Proof:} According to Theorem \ref{PT}, the law of $V(B^{\lambda,br})$ is uniquely determined by that of the triple $(Z^{\lambda},B^{ex,Z^{\lambda}},F^{\lambda,1-Z^{\lambda}})$. It suffices to prove that the law of the process  on the left hand side of \eqref{5} is entirely characterized by the same triple. Following Theorem $3.1$ in Perman et al \cite{PPY} and the discussion below Lemma $4.10$ of Pitman \cite{Pitman}, conditioned on $\Delta:=D_U-G_U$, $(\tilde{X}_t; 0 \leq t \leq \Delta)$ and $(\tilde{X}_t; \Delta \leq t \leq 1)$ are independent and $\Delta$ corresponds to the length of first excursion of  $X$ via length-biased sampling:
$f_{\Delta}(t)=f_{Z^{\lambda}}(t)~\mbox{as in \eqref{aldous}}.$
Thus, $\Delta \stackrel{d}{=} Z^{\lambda}$. Finally, since $L_t^0(\tilde{X})=0$ on $(0,D_U-G_U)$, we have that $(|\tilde{X}(t)|-L_t^0(\tilde{X}); 0 \leq t \leq D_U-G_U |L_1^0(\tilde{X})=|\lambda|)$ is a Brownian excursion of length $\Delta$ (conditionally) independent of $(|\tilde{X}(t)|-L_t^0(\tilde{X});  D_U-G_U \leq t \leq 1 |L_1^0(\tilde{X})=|\lambda|)$, which is a first passage bridge through level $\lambda$ of length $1-\Delta$ by construction \eqref{4}.  $\square$\\\\
\textbf{Remark:} The process $X$ defined in the above theorem is a Brownian bridge conditioned on its local time, see Chassaing and Janson \cite{CJ} for detail discussions. In addition, the proof of Theorem $3.1$ in Perman et al \cite{PPY} is extensively based on the concept of Palm distribution, which can be read from Fitzsimmons et al \cite{FPY}. 
\subsection{Vervaat bridges are not Markov}
It is natural to ask whether the Vervaat bridges are Markov (with respect to their induced filtrations). In the case of negative endpoints, it is equivalent to ask whether the entrance law after the excursion piece is nice enough for the first passage bridge to produce Markov property. The following result gives a negative answer.
\begin{proposition}
\label{notmarkov}
$(V(B^{\lambda,br})_t; 0 \leq t \leq 1))$ where $\lambda<0$ is not Markov with respect to its induced filtration.
\end{proposition}

Before proving the proposition, we introduce some notations that we use in the current section and rest of the paper.
For $x,y >0$, denote
\begin{displaymath}
\tilde{q}_{t}(x,y):=\frac{1}{xy\sqrt{2 \pi t}} \left(\exp \left(-\frac{(x-y)^2}{2t}\right)- \exp \left(-\frac{(x+y)^2}{2t}\right)\right).
\end{displaymath}
Note that $\tilde{q}_{t}(x,y) y^{2}\,dy$ is the transition kernel of three dimensional Bessel process and
\begin{displaymath}
\tilde{q}_{t}(0,y)=\lim_{x\rightarrow 0^{+}}\tilde{q}_{t}(x,y)=
\dfrac{2}{\sqrt{2\pi t^{3}}}\exp\left(-\dfrac{y^{2}}{2t}\right)=\dfrac{2}{y}g_t(y),
\end{displaymath}
\begin{displaymath}
\tilde{q}_{t}(0,0)=
\dfrac{2}{\sqrt{2\pi t^{3}}},
\end{displaymath}
where $g_t(y)$ is the density of the first hitting at level $y$ for Brownian motion given in \eqref{Eq4}.\\\\
\textbf{Proof of Proposition \ref{notmarkov}:} Fix $t_0 \in (0,1)$ and $x_0>0$. Let $T_{t_0}$ be the first return of $V(B^{\lambda,br})$ to $0$ after time $t_0$. Consider the distribution of $T_{t_0}$ given $V(B^{\lambda,br})_{\frac{t_0}{2}}=0$ and $V(B^{\lambda,br})_{t_0}=x_0$. According to Theorem \ref{PT}, given $T_{t_0}$, $(V(B^{\lambda,br})_t; t_0 \leq t \leq T_{t_0})$ and $(V(B^{\lambda,br})_t; T_{t_0} \leq t \leq 1)$ are two independent first passage bridges from $x_0>0$ to $0$ respectively from $0$ to $\lambda<0$. Therefore, its density is given by
\begin{align}
\label{new1}
f_1(t)&=\frac{g_{t-t_0}(x_0) g_{1-t}(|\lambda|)}{g_{1-t_0}(1+|\lambda|)}1_{t > t_0} \notag\\
         &=\frac{C_1(t_0,x_0,\lambda)}{\sqrt{(t-t_0)^3(1-t)^3}} \exp \left( -\frac{x_0^2}{2(t-t_0)}-\frac{\lambda^2}{2(1-t)} \right) 1_{t>t_0},
\end{align}
for some $C_1(t_0,x_0,\lambda)>0$. Next we consider the distribution of $T_{t_0}$ given that $\forall u \in (0,t_0), V(B^{\lambda,br})_u>0$ and $V(B^{\lambda,br})_{t_0}=x_0$, whose density can be computed using Bayes recipe:
\begin{align}
\label{new2}
f_2(t)&=C_2(t_0,x_0,\lambda) \frac{\tilde{q}_{t_0}(0,x_0) \tilde{q}_{t-t_0}(x_0,0)}{\tilde{q}_t(0,0)}x_0^2 f_{Z^{\lambda}}(t)1_{t>t_0} \notag\\
         &=C_2^{'}(t_0,x_0,\lambda) \frac{t}{\sqrt{(t-t_0)^3(1-t)^3}}\exp \left( -\frac{x_0^2}{2(t-t_0)}-\frac{\lambda^2}{2(1-t)} \right) 1_{t>t_0}.
\end{align}
for some $C_2(t_0,x_0,\lambda)>0$ and $C_2^{'}(t_0,x_0,\lambda)>0$.
Comparing \eqref{new1} to \eqref{new2}, we have that $f_2(t)=C_{1,2}(t_0,x_0,\lambda)tf_1(t)$ for some $C_{1,2}(t_0,x_0,\lambda)>0$. The two conditional densities of $T_{t_0}$ fail to be equal and we conclude that the Vervaat bridges with negative endpoint are not Markov.  $\square$
\begin{center}
\includegraphics[width=0.6\textwidth]{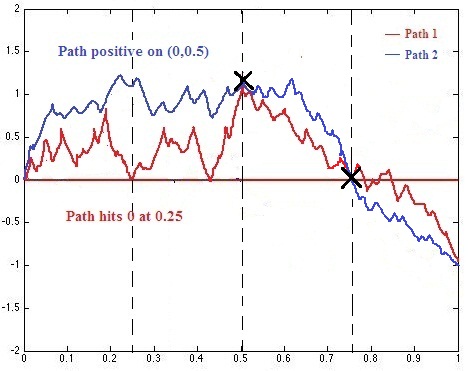}\\
Fig 5. Paths which make the Vervaat bridge non-Markov, $t_0=0.5$ and $x_0=1$.
\end{center}
\textbf{Remark:} The counter-example provided in the proof of Proposition \ref{notmarkov} indicates that the main reason that makes Vervaat bridges with negative endpoint non-Markov is the lack of information on $Z$. Indeed, for $s \leq t \leq 1$, $V(B^{\lambda,br})_t$ depends not only on $V(B^{\lambda,br})_s$ but also on the event $\{Z \leq s\}$. 

It is well-known that the time reversal of any Markov process is still Markov. This result leads to the following corollary saying that the Vervaat bridges with positive endpoint is not Markov as well.
\begin{corollary}
$(V(B^{\lambda,br})_t; 0 \leq t \leq 1))$ where $\lambda>0$ is not Markov with respect to its induced filtration.
\end{corollary}

Now we know that the Vervaat bridges of non-zero endpoint are not Markov. Thus it is natural to ask how bad they may behave so that Markov property cannot be produced. This leads to the question that whether they are semimartingales. The discussion of this question is reported later to the Section $4.2$.  
\subsection{Relation with drifting excursion}
In Bertoin \cite{Bertoin}, he studied a fragmentation process by considering the excursion dragged down by drift $\lambda<0$:
$$B_t^{ex, \lambda \downarrow}:=B^{ex}_t+\lambda t, ~\mbox{for}~0 \leq t \leq 1.$$

Notice that $V(B^{\lambda,br})$ (with $\lambda<0$) also looks similar to this process except that the former always stays above the line $t \rightarrow \lambda t$ while the latter doesn't share this property. A natural way to relate these two processes is to see whether conditioned on staying above the dragging line, the Vervaat bridge is absolutely continuous with respect to drifting excursion. First we need to justify that the conditioning event has positive probability. The next proposition provides a positive answer with an explicit formula.
\begin{proposition}
\label{comp}
$\forall \lambda<0$,
$$ \mathbb{P}(\forall t \in (0,1),~V(B^{\lambda,br})_t >\lambda t)=1-|\lambda| \exp\left(\frac{\lambda^2}{2}\right) \int_{|\lambda|}^{\infty} \exp\left(-\frac{t^2}{2}\right)dt.$$
\end{proposition}
\textbf{Proof:} Following from Proposition $15$ of Schweinsberg \cite{Schweinsberg}, fix $x \in [\lambda,0]$ we know the probability of a first passage bridge through level $\lambda$ to stay above the dragging line tying $x$ to $\lambda$:
\begin{equation}
\label{15}
\mathbb{P}(\forall t \in [0,l], F^{\lambda,l}(t)>x-(x-\lambda)t)=\frac{|x|}{|\lambda|}.
\end{equation}
Therefore,
\begin{align}
\mathbb{P}(\forall t \in (0,1),~V(B^{\lambda,br})_t >\lambda t) &=\int_0^1  \mathbb{P}\left(\forall s \in (t,1),~V(B^{\lambda,br})_s >\lambda s | Z^{\lambda}=t\right)f_{Z^{\lambda}}(t)dt  \notag\\
                                                                                                    &=\int_0^1  t \frac{|\lambda|}{\sqrt{2 \pi t(1-t)^3}} \exp\left(-\frac{\lambda^2 t}{2(1-t)}\right)dt   \notag\\
                                                                                                    &=\mathbb{E}Z^{\lambda}.   \notag
\end{align}
where the first equality follows from the fact that the excursion piece is always above the dragging line and the second equality is a direct consequence of \eqref{15}. Following the notations of discussion below Lemma $4.10$ in Pitman \cite{Pitman},
$$\mathbb{E}Z= h_{-2}(\lambda) \mathbb{E}B_1^2=1-\lambda \exp\left(\frac{\lambda^2}{2}\right) \int_{\lambda}^{\infty} \exp\left(-\frac{t^2}{2}\right)dt,$$
where $h_{-2}$ is the Hermite function of index $-2$.  $\square$

Now we know that the Vervaat bridge (with negative endpoint) conditioned to stay above the dragging line is well-defined. In addition, the law of its first return to $0$ is given by:
\begin{equation}
\label{16}
f_{\widetilde{Z}^{\lambda}}(t)=\frac{t}{1-|\lambda| \exp\left(\frac{\lambda^2}{2}\right) \int_{|\lambda|}^{\infty} \exp\left(-\frac{t^2}{2}\right)dt}f_{Z^{\lambda}}(t). 
\end{equation}

The next theorem provides a path decomposition result of Vervaat's bridge conditioned to stay above the dragging line and establishes connection to drifting excursion.
\begin{theorem}
Let $\lambda<0$. Given $\widetilde{Z}^{\lambda}$ the length of first excursion of $(V(B^{\lambda,br})_t; 0 \leq t \leq 1|\forall t \in (0,1), V(B^{\lambda,br})_t>\lambda t)$ (whose distribution density is given by \eqref{16}), the path is decomposed into two (conditionally) independent pieces:
\begin{itemize}
\item
$\left(V(B^{\lambda,br})_u;0 \leq u \leq \widetilde{Z}^{\lambda}|\forall t \in (0,1), V(B^{\lambda,br})_t>\lambda t\right)$ is an excursion of length $\widetilde{Z}^{\lambda}$;
\item
 $\left(V(B^{\lambda,br})_u;\widetilde{Z}^{\lambda} \leq u \leq 1|\forall t \in (0,1), V(B^{\lambda,br})_t>\lambda t\right)$ is a first passage bridge of length $1-\widetilde{Z}^{\lambda}$ conditioned to stay above $t \rightarrow \lambda(t+\widetilde{Z}^{\lambda})$ for $t \in (0,1-\widetilde{Z}^{\lambda})$.
\end{itemize}
In addition, $(V(B^{\lambda,br})_t; 0 \leq t \leq 1|\forall t \in (0,1), V(B^{\lambda,br})_t>\lambda t)$ is absolutely continuous with respect to $(B_t^{ex,\lambda \downarrow}; 0\leq t \leq 1)$. The corresponding density is: $$\frac{H}{1-|\lambda| \exp\left(\frac{\lambda^2}{2}\right) \int_{|\lambda|}^{\infty} \exp\left(-\frac{t^2}{2}\right)dt},$$
where $H:=\inf \{ t>0; B^{ex, \lambda \downarrow}_{t}<0\}$.
\end{theorem}
\textbf{Proof:} According to Proposition $11$ of Bertoin \cite{Bertoin}, $H$ is distributed as \eqref{aldous}. Following Theorem $2.6$ of Chassaing and Jason \cite{CJ}, conditioned on $H$, $(B^{ex, \lambda \downarrow}_t; 0 \leq t \leq H)$ is a Brownian excursion of length $H$. In addition, Proposition $4$ of Schweinsberg \cite{Schweinsberg} states that given $H$, $(B^{ex, \lambda \downarrow}_t; H \leq t \leq 1)$ is a first passage bridge of length $1-H$ conditioned to stay above the line $t \rightarrow \lambda(t+H)$ for $t \in (0,1-H)$, (conditionally) independent of the excursion piece. By change of measures, we obtain the same triple characterization in law. $\square$
\subsection{Convex minorant of Vervaat bridges}
In this part, we will study some properties of convex minorant of Vervaat bridge $V(B^{\lambda,br})$ where $\lambda<0$. The convex minorant of a real-valued function $(X_t; t \in [0,1])$ is the maximal convex function $(C_t; t \in [0,1])$ such that $ \forall t \in [0,1], C_t \leq X_t $. We refer to the points where the convex minorant equals  the process as vertices. Note that these points are also the endpoints of the linear segments. See Pitman and Ross \cite{PR} and Abramson et al \cite{Abr} for general background.

Similar to the computation in Proposition \ref{comp} , we have the explicit formula for the distribution of the last segment's slopes.
\begin{corollary}
Denote $s_l$ the slope of the last segment of the convex minorant for $(V(B^{\lambda,br})_t;$ $0 \leq t \leq 1)$. $\forall a \in [\lambda,0]$, we have
$$\mathbb{P}(s_l \in [\lambda, a])=1+a \exp\left(\frac{\lambda^2}{2}\right) \int_{|\lambda|}^{\infty} \exp\left(-\frac{t^2}{2}\right)dt.$$
\end{corollary}

As discussed in Pitman and Ross \cite{PR}, a standard first passage bridge can only have accumulations of linear segments at its start point (while Brownian motion has accumulations at two endpoints). However, seen in the beginning of the section, the greatest difference between the Vervaat bridges and the standard first passage bridges is the first excursion piece for the former. Then we can expect that the Vervaat bridges have almost surely a finite number of segments.
\begin{proposition}
The number of segments of the convex minorant of  $V(B^{\lambda,br})$ for $\lambda<0$ is a.s. finite.
\end{proposition}
\textbf{Proof:} We adopt a sample paths argument. Consider a sample path of Brownian bridge $B^{\lambda,br}$ where $\lambda<0$and $1-A^{\lambda}:=\argmin B^{\lambda,br}$ (which is a.s. unique). Note that $V(B^{\lambda,br})_t >0$ for $t \in (0,A]$. Consequently, the first vertex of the Vervaat bridge $\alpha_1>A$ a.s. According to Pitman and Ross \cite{PR}, there can be only a finite number of segments on $[\alpha_1,1]$ since accumulations can only happen at $0$ on the restricted path $B^{\lambda,br}|_{[0,1-A]}$. Thus, the number of segments of the Vervaat bridges is a.s. finite.  $\square$

However, we expect a stronger result regarding the number of segments:
\begin{conjecture}
The expected number of segments of the Vervvat bridges is finite.
\end{conjecture}
\section{The Vervaat transform of Brownian motion}
In this section, we study the Vervaat transform of Brownian motion. We first prove that the process is not Markov with respect to its induced filtration. Next, it is shown to be a semimartingale. Finally, we provide the mean and the variance of this process.
\subsection{$V(B)$ is not Markov}
An important property of $V(B)$ is that it has the same terminal value as $B$: $V(B)_1=B_1$. We have two cases: If $B_1>0$, then $V(B)$ never returns to $0$ along the path. Otherwise $B_1 \leq 0$, then $V(B)_1 \leq 0$. By path continuity, $V(B)$ has to hit $0$ somewhere on its path.
\begin{proposition}
$(V(B)_t; 0 \leq t \leq 1)$ is not Markov with respect to its induced filtration.
\end{proposition}
\textbf{Proof:} According to the above discussion, 
\begin{equation}
\label{17}
\mathbb{P}(V(B)_1>0|V(B)_{\frac{1}{2}}>0, V(B)_{\frac{1}{4}}=0)=0;
\end{equation}
since once it hits $0$ on its path, $V(B)$ has to end negatively. On the other hand,
\begin{align}
\label{18}
\mathbb{P}(V(B)_1>0|\forall t \in (0, \frac{1}{2}],V(B)_t>0)&=\frac{\mathbb{P}(V(B)_1>0~\mbox{and}~\forall t \in (0, \frac{1}{2}],V(B)_t>0)}{\mathbb{P}(\forall t \in (0, \frac{1}{2}],V(B)_t>0)} \notag\\
                                                                                            &\geq \mathbb{P}(\forall t \in([0,1],V(B)_t>0) \notag\\
                                                                                            &=\mathbb{P}(B_1>0)=\frac{1}{2}. 
\end{align}
\begin{center}
\includegraphics[width=0.6\textwidth]{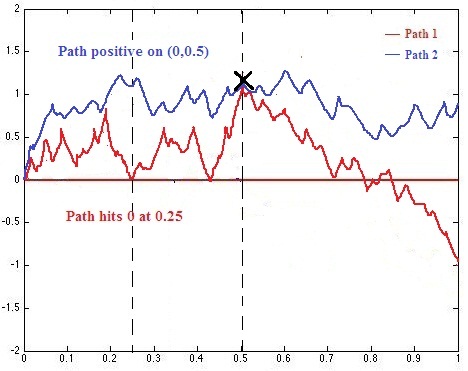}\\
Fig 6. Paths make Vervaat's transform of BM non-Markov.
\end{center}
By comparing \eqref{17} and \eqref{18}, we see that these two conditional probabilities fail to be equal, which implies that $(V(B)_t; 0 \leq t \leq 1)$ is not Markov.   $\square$ \\\\
\textbf{Remark:} If we denote $T:=\inf\{t>0; V(B)_t=0\}$, we have $\{T \leq 1\}=\{V(B)_1 \leq 0\}$. Formally this means that we obtain the information at time $1$ from some prior time, which violates the Markov property.
\subsection{$V(B)$ is a semimartingale}
In general, when a process is Markov (with state space in $\mathbb{R}^d$), we know sufficient and necessary conditions for it to be a semimartingale, see Cinlar et al \cite{CJSP}. However, we have seen in the preceding subsection that $V(B)$ is not Markov. Therefore, whether $V(B)$ is a semimartingale or not cannot be judged by classical Markov-semimartingale procedures. In this section, we provide a soft argument to prove that $V(B)$ is indeed a semimartingale with respect to its induced filtration using Denisov's decomposition for Brownian motion as well as Bichteler-Dellacherie's characterization for semimartingales.

We first recall some paths decomposition result for standard Brownian motion, which permits a characterization for the Vervaat transform. Following the notations in the introduction, $A$ is the a.s. arcsine split $(1-A:=\argmin_{t \in [0,1]} B_t)$ for a standard Brownian motion. The following theorem is due to Denisov \cite{Denisov}:
\begin{theorem} \label{Denisov} \textbf{Denisov's decomposition} \cite{Denisov}
 Given $A$ (which is arcsine distributed, i.e. $f_A(a)=\frac{1}{\pi \sqrt{a(1-a)}}$), the path is decomposed into two independent pieces:
\begin{itemize}
\item
$\left(\frac{1}{\sqrt{A}}(B_{1-A+uA}-B_{1-A});0 \leq u \leq 1\right)$  is a standard Brownian meander;
\item
$\left(\frac{1}{\sqrt{1-A}}(B_{(1-A)(1-u)}-B_{1-A});0 \leq u \leq 1\right)$ is a standard Brownian meander.
\end{itemize}
\end{theorem}
\textbf{Remark:} The theorem simply says that given its a.s. minimum $A$, a standard Brownian motion is split into two conditional independent meanders of length $A$ and $1-A$ joint back to back. Therefore, $(V(B)_t; 0 \leq t \leq 1)$ can be viewed as mixing of two independent joint back-to-back Brownian meanders with respect to arcsine distribution.

Now we turn to some results in the classical semimartingale theory.
Given a filtration $(\mathcal{F}_t)_t$, a process $H$ is said to be simple predictable if $H$ has a representation
$$\forall t \in [0,1],~H_t=H_01_{\{0\}}(t)+\sum_{i=1}^{n-1} H_i 1_{(t_i,t_{i+1}]}(t)$$
where $H_i \in \mathcal{F}_{t_i}$ and $|H_i|< \infty$ a.s. for $0=t_1 \leq t_2 \leq...\leq t_n \leq \infty$. \\

Denote $\mathcal{S}$ the collection of simple predictable processes and $\mathcal{B}=\{H \in \mathcal{S}: |H| \leq 1\}$ (the unit ball in $\mathcal{S}$). For a given process, we define a linear mapping $I_X: \mathcal{S} \rightarrow \mathbb{L}^0$ by
$$I_X(H_t)=H_0X_0+\sum_{i=1}^{n-1}H_i(X_{t_i}-X_{t_{i-1}}),$$
for $H\in \mathcal{S}$. In fact, $I_X$ is defined as stochastic integral with respect to $X$ for simple predictable processes.

The following theorem, proved independently by Bichteler \cite{Bich} and Dellacherie \cite{Dellacherie} provides a useful characterization for semimartingales. We refer the readers to Jacod \cite{Jacod}, Protter \cite{Protter} and Rogers and Williams \cite{RW} for more details.
\begin{theorem}
\label{BD}
\textbf{Bichteler-Dellacherie's theorem} \cite{Bich}, \cite{Dellacherie} An adapted, c\`{a}dl\`{a}g process $X$ is a semimartingale iff $I_X(\mathcal{B})$ is bounded in probability, which means
$$\lim_{\eta \rightarrow \infty} \sup_{H \in \mathcal{B}} \mathbb{P}(|I_X(H)| \geq \eta)=0.$$
\end{theorem}
\textbf{Remark:} Fundamentally, this theorem tells that the notion of semimartingale is equivalent to the notion of "good stochastic integrator" and it depends only on the law of the processes.

We now state the main theorem of the section:
\begin{theorem}
\label{semi}
$(V(B)_t; 0 \leq t \leq 1)$ is semimartingale with respect to its induced filtration.
\end{theorem}
\textbf{Proof:} Fix $H \in \mathcal{B}$ and $\eta>0$,
\begin{equation}
\label{19}
\mathbb{P}(|I_{V(B})(H)|>\eta)=\int_0^1 \mathbb{P}(|I_{V(B)|A=a}(H)|>\eta) \frac{1}{\pi \sqrt{a(1-a)}}da.
\end{equation}
Note that $(V(B)|A=1)$ is a standard Brownian meander and $B^{me}\stackrel{d}{=} R^{0 \rightarrow \rho}$ where $R^{x \rightarrow y}$ is a three dimensional Bessel bridge from $x$ to $y$ and $\rho$ is Rayleigh distributed $\mathbb{P}(\rho \in dx)=x \exp(-\frac{x^2}{2})dx$. Girsanov's change of measure theorem guarantees that $(V(B)|A=1)$ is a semimartingale and so is $(V(B)|A=0)$ (see e.g. Imhof \cite{Imhof} and Az\'{e}ma-Yor \cite{AY}). Thus, by Theorem \ref{BD}, 
$$\lim_{\eta \rightarrow \infty} \sup_{H \in \mathcal{B}} \mathbb{P}(|I_{V(B)|A=1}(H)| \geq \eta)=0;$$
$$\lim_{\eta \rightarrow \infty} \sup_{H \in \mathcal{B}} \mathbb{P}(|I_{V(B)|A=0}(H)| \geq \eta)=0.$$
From \eqref{19}, to prove $\sup_{H \in \mathcal{B}} \mathbb{P}(|I_{V(B})(H)| \geq \eta) \rightarrow 0$, we need some uniform control for $\sup_{H \in \mathcal{B}} \mathbb{P}(|I_{V(B)|A=a}(H)|>\eta)$. The following result permits to estimate for all $a \in [0,1]$, $\sup_{H \in \mathcal{B}} \mathbb{P}|I_{V(B)|A=a}(H)|>\eta)$ in terms of  $\sup_{H \in \mathcal{B}} \mathbb{P}(|I_{V(B)|A=1}(H)|>\eta)$ and $\sup_{H \in \mathcal{B}} \mathbb{P}(|I_{V(B)|A=0}(H)|>\eta)$.
\begin{lemma}
\label{control}
$\forall a \in (0,1]$,
$$\sup_{H \in \mathcal{B}} \mathbb{P}(|I_{V(B)|A=a}(H)| > \eta) \leq  \sup_{H \in \mathcal{B}} \mathbb{P}(|I_{V(B)|A=1}(H)| > \frac{\eta}{2})+\sup_{H \in \mathcal{B}} \mathbb{P}(|I_{V(B)|A=0}(H)| > \frac{\eta}{2}).$$
\end{lemma}
\textbf{Proof:} Observe that $I_{V(B)|A=a}(H)=I_1+I_2$, where
$$I_1:=\sum_i H_i(V(B)_{\tau_{i+1} \wedge a}-V(B)_{\tau_{i+1} \wedge a}) \quad \mbox{and} \quad I_2:=\sum_i H_i(V(B)_{\tau_{i+1} \vee a}-V(B)_{\tau_{i+1} \vee a}).$$
We have then,
$$\mathbb{P}(|I_{V(B)|A=a}(H)|>\eta) \leq \mathbb{P}(|I_1| >\frac{\eta}{2})+\mathbb{P}(|I_2| >\frac{\eta}{2}).$$
Denote $\tilde{I_1}=\frac{I_1}{\sqrt{a}}$ and note that 
\begin{align}
\tilde{I_1}&=\sum_i H_i \frac{V(B)_{\tau_{i+1} \wedge a} -V(B)_{\tau_i \wedge a}}{\sqrt{a}} \notag\\
                 &=\sum_i H_i (\tilde{V}(B)_{\frac{\tau_{i+1} \wedge a}{a}}-\tilde{V}(B)_{\frac{\tau_i \wedge a}{a}}). \notag
\end{align}
where $(\tilde{V}(B)_t; 0 \leq t \leq 1)$ is a standard Brownian meander and $H_i$ is $\mathcal{F}^{\tilde{V}}_{\frac{\tau_i \wedge a}{a}}$-adapted $\forall i$.\\
We have then 
$$\mathbb{P}(|I_1| > \frac{\eta}{2}) \leq \mathbb{P}(|\tilde{I_1}| >\frac{\eta}{2}) \leq \sup_{H \in \mathcal{B}} \mathbb{P}(|I_{V(B)|A=1}(H)| > \frac{\eta}{2}).$$
Similarly, by independence of two decomposed meanders,
$$\mathbb{P}(|I_2| > \frac{\eta}{2}) \leq \mathbb{P}(|\tilde{I}_2| > \frac{\eta}{2})=\sup_{H \in \mathcal{B}} \mathbb{P}(|I_{V(B)|A=0}(H)| > \frac{\eta}{2}).$$ 
where $\tilde{I}_2$ is the stochastic integral associated to reversed Brownian meander.  $\square$\\\\
Now return to the proof of Theorem \ref{semi}. According to Lemma \ref{control},
\begin{align}
&~~~~\sup_{H \in \mathcal{B}} \mathbb{P}(|I_X(H)|>\eta) \notag\\
&=\sup_{H \in \mathcal{B}}\int_0^1 \mathbb{P}(|I_{V(B)|A=a}(H)|>\eta) \frac{1}{\pi \sqrt{a(1-a)}}da \notag\\
& \leq  \int_0^1 [\sup_{H \in \mathcal{B}} \mathbb{P}(|I_{V(B)|A=1}(H)|>\frac{\eta}{2})+\sup_{H \in \mathcal{B}} \mathbb{P}(|I_{V(B)|A=0}(H)|>\frac{\eta}{2})]\frac{1}{\pi \sqrt{a(1-a)}}da \notag
\end{align}
which goes to $0$ as $\eta \rightarrow \infty$.  $\square$

The following corollary states that the Vervaat bridges are also semimartingales, which provides an alternative proof of the semimartingale property for Vervaat bridges obtained in Section $3.3$.
\begin{corollary}
For each fixed $\lambda \in \mathbb{R}$, $(V(B^{\lambda,br})_t; 0 \leq t \leq 1)$ is a semimartingale with respect to its induced filtration.
\end{corollary}
\textbf{Proof:} Fix $H \in \mathcal{B}$ and $\eta>0$,
\begin{equation}
\label{20}
\mathbb{P}(I_{V(B)} (H)>\eta)=\int_{\mathbb{R}}\mathbb{P}(I_{V(B^{\lambda,br})}(H) > \eta)\frac{1}{\sqrt{2 \pi}} \exp\left(-\frac{\lambda^2}{2}\right)d\lambda.
\end{equation}
Note that $V(B^{0,br})$ is Brownian excursion, thus a semimartingale. It suffices then to prove $V(B^{\lambda,br})$ for $\lambda \neq 0$ is a semimartingale. If not the case, $\exists \epsilon>0$, $\forall K>0~\exists \eta>K$ such that 
$$\sup_{H \in \mathcal{B}}\mathbb{P}(I_{V(B^{\lambda,br})}(H) > \eta)>\epsilon.$$
Note in addition that $(H, \lambda) \rightarrow \mathbb{P}(I_{V(B^{\lambda,br})}(H) > \eta)$ is jointly continuous in $\mathcal{B} \times (\mathbb{R}\setminus \{0\})$ (concatenation of continuity that is left for readers  to check). Thus $\exists H_{\lambda,\epsilon} \in \mathcal{B}$ and $\theta \in (0,|\lambda|)$ such that $\forall \tilde{\lambda} \in (\lambda-\theta,\lambda+\theta)$,
\begin{equation}
\label{21}
\mathbb{P}(I_{V(B^{\tilde{\lambda},br})}(H) > \eta)>\frac{\epsilon}{2}.
\end{equation}
Injecting \eqref{21} into \eqref{20}, we obtain:
$$\mathbb{P}(I_{V(B)}(H) > \eta)>\frac{\epsilon}{2}\int_{\lambda-\theta}^{\lambda+\theta}\frac{1}{\sqrt{2 \pi}} \exp\left(-\frac{\lambda^2}{2}\right)d\lambda.$$
which violates that fact that $V(B)$ is a semimartingale.  $\square$

However, one can hardly derive an explicit decomposition formula using Bichteler-Dellacherie's approach. Let's explain why: a generic approach for the proof of Bichteler-Dellacherie's theorem is to find $\mathbb{Q}$ equivalent to $\mathbb{P}$ such that $X$ is $\mathbb{Q}-$quasimartingale (see e.g. Protter \cite{Protter} for definition). By Rao's theorem, $X$ is $\mathbb{Q}-$semimartingale, which is also $\mathbb{P}-$semimartingale by Girsanov's theorem. Note that Rao's theorem is based on Doob-Meyer's decomposition theorem, which in general does not give an explicit expression for two decomposed terms (in fact they are defined as some limiting processes). 
\subsection{Expectation and variance for $V(B)$}
In the current subsection, we provide the formulae for the first two moments of the Vervaat transform of Brownian motion.
\begin{proposition}
$\forall t \in [0,1]$, we have:
\begin{equation}
\label{23}
\mathbb{E}V(B)_t=\sqrt{\frac{8}{\pi}} (\sqrt{t}+\sqrt{1-t}-1);
\end{equation}
\begin{equation}
\label{24}
\mathbb{E}(V(B)_t^2)=3t+\frac{4-8t}{\pi} \arcsin\sqrt{t}-\frac{4}{\pi}\sqrt{t(1-t)}.
\end{equation}
\end{proposition}
\begin{center}
\includegraphics[width=0.6\textwidth]{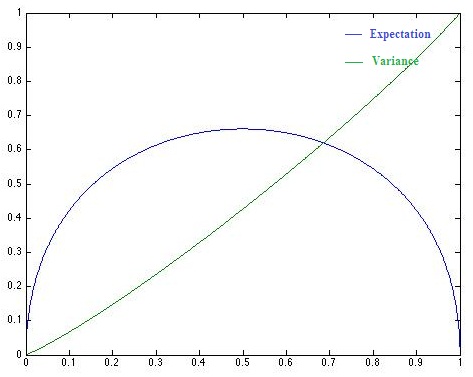}\\
Fig 7. Expectation and variance for $V(B)$.
\end{center}

The computation is based on Theorem \ref{Denisov}, which is stated in the Section $4.2$ as well as the following identities for standard Brownian meander, whose proof will be reported to the Appendix:
\begin{proposition}
\label{compmean}
Let $(B^{me}_t, t \in [0,1])$ be standard Brownian meander. We have:
\begin{equation}
\label{25}
\mathbb{E}B^{me}_t=\sqrt{\frac{2}{\pi}} (\sqrt{t(1-t)}+ \arcsin \sqrt{t}).
\end{equation}
\begin{equation}
\label{26}
\mathbb{E}(B^{me}_t)^2=3t-t^2.
\end{equation}
\begin{equation}
\label{27}
\mathbb{E}B^{me}_tB^{me}=2 \sqrt{t}.
\end{equation}
\end{proposition}
\textbf{Remark:} One can also think of computing the expectation and the variance of the Vervaat bridges. However, we are not able to derive some explicit formulae for them except in the case of zero endpoint (correspond to Brownian excursion). Also note that the expectation as well as the variance of the Vervaat transform of Brownian motion can be obtained by discrete approximation.
\subsubsection{Expectation for $V(B)$}
Let's compute the expectation of $(V(B)_t, 0 \leq t \leq 1)$. Recall that $A$ is the a.s. arcsine split $(1-A:=\argmin_{t \in [0,1]} B_t)$, we have:
$$\mathbb{E}V(B)_t=\mathbb{E}(V(B)_t 1_{A>t})+ \mathbb{E}(V(B)_t 1_{A \leq t}).$$
\begin{lemma}
\label{pre1}
\begin{equation}
\label{28}
\forall t \in [0,1],~\mathbb{E}(V(B)_t 1_{A>t})=\sqrt{\frac{2}{\pi}}(\sqrt{1-t}+2 \sqrt{t}-t-1).
\end{equation}
\end{lemma}
\textbf{Proof:} According to the formula for $\mathbb{E}B^{me}(t), t\in [0,1]$, we have:
\begin{align}
\mathbb{E}(V(B)_t 1_{A>t})&=\int_t^1 \frac{\sqrt{a} \mathbb{E}B^{me}(\frac{t}{a})}{\pi \sqrt{a(1-a)}} da\notag\\
                                         &\stackrel{\eqref{25}}{=}\frac{\sqrt{2}}{\pi^{\frac{3}{2}}} \int_t^1 \frac{\sqrt{\frac{t}{a}} \sqrt{1-\frac{t}{a}} + \arcsin \sqrt{\frac{t}{a}}}{\sqrt{1-a}}   \notag\\
                                         &~=\frac{\sqrt{2}}{\pi^{\frac{3}{2}}} (\alpha_1+\alpha_2),  \notag                               
\end{align} 
where $\alpha_1:=\sqrt{t} \int_t^1 \frac{\sqrt{a-t}}{a \sqrt{1-a}}da$ and $\alpha_2:=\int_t^1 \frac{\arcsin\sqrt{\frac{t}{a}}}{\sqrt{1-a}}da.$
Using integration by parts, we get:
\begin{align}
\alpha_2&= \left[-2 \sqrt{1-a} \arcsin\sqrt{\frac{t}{a}}\right]_t^1 + 2 \int_t^1 \sqrt{1-a} \left(\arcsin\sqrt{\frac{t}{a}}\right)^{'} da \notag\\
     &=\pi \sqrt{1-t} -\sqrt{t} \int_t^1 \frac{\sqrt{1-a}}{a \sqrt{a-t}}da.   \notag                                                                                       
\end{align} 
Therefore,
\begin{align}
\label{29}
\alpha_1+\alpha_2&= \pi \sqrt{1-t} + \sqrt{t} \int_t^1  \left(\frac{\sqrt{a-t}}{a \sqrt{1-a}} - \frac{\sqrt{1-a}}{a \sqrt{a-t}}\right) da \notag\\
             &=\pi \sqrt{1-t} + 2 \sqrt{t} \int_t^1 \frac{da}{\sqrt{1-a}\sqrt{a-t}} - \sqrt{t}(t+1) \int_t^1 \frac{da}{a \sqrt{1-a}\sqrt{a-t}}.
\end{align} 
By change of variables, we obtain:
\begin{equation}
\label{30}
\int_t^1 \frac{da}{\sqrt{1-a}\sqrt{a-t}}\stackrel{a=t+(1-t)x}{=}\int_0^1 \frac{dx}{\sqrt{x(1-x)}}=\pi.
\end{equation}
\begin{align}
\label{31}
\int_t^1 \frac{da}{a \sqrt{1-a}\sqrt{a-t}} &\stackrel{a=t+(1-t)x}{=} \int_0^1 \frac{dx}{\sqrt{x(1-x)} [t+(1-t)x]} \notag\\
                                                                     &~~\stackrel{x=\sin^2 \theta}{=} \frac{2}{1-t}\int_0^{\frac{\pi}{2}} \frac{d \theta}{\sin^2 \theta + \frac{t}{1-t}}.
\end{align}
Observe that $\theta \rightarrow \frac{1}{\sqrt{a(a+1)}} \arctan (\sqrt{\frac{a}{a+1}}\tan x)$ is a primitive $\theta \rightarrow \frac{1}{\sin^2 \theta +a}$ for $a>0$. Take $a=\frac{t}{1-t}$, we have from \eqref{31}:
$$\forall t>0,~\int_t^1 \frac{da}{a \sqrt{1-a}\sqrt{a-t}}=\frac{\pi}{\sqrt{t}}.$$
and consequently,
\begin{equation}
\label{32}
\forall t \geq 0,~\sqrt{t} \int_t^1 \frac{da}{a \sqrt{1-a}\sqrt{a-t}}=\pi 
\end{equation}
Combining \eqref{29}, \eqref{30} and \eqref{32}, we get \eqref{28}.  $\square$

\begin{lemma}
\label{pre2}
\begin{equation}
\label{33}
\forall t \in [0,1],~\mathbb{E}(V(B)_t 1_{A \leq t})=\sqrt{\frac{2}{\pi}} (\sqrt{1-t}+t-1).
\end{equation}
\end{lemma}
\textbf{Proof:} We have $\mathbb{E}(V(B)_t 1_{A \leq t})=\mathbb{E}((V(B)_t-V(B)_1)1_{A \leq t})+\mathbb{E}(V(B)_1 1_{A \leq t})$, where the first term can be derived from \eqref{28} by change of variables:
\begin{equation}
\label{34}
\mathbb{E}((V(B)_t-V(B)_1)1_{A \leq t})=\sqrt{\frac{2}{\pi}}(\sqrt{t}+2 \sqrt{1-t}+t-2).
\end{equation}
The second one can be computed by independence of two meanders:
\begin{align}
\label{35}
\mathbb{E}(V(B)_1 1_{A \leq t}) &= \int_0^t \mathbb{E}(V(B)_1|A=a) \frac{da}{\pi \sqrt{a(1-a)}} \notag\\
                                                    &=\int_0^t \left(\sqrt{\frac{\pi}{2}a}-\sqrt{\frac{\pi}{2}(1-a)}\right)\frac{da}{\pi \sqrt{a(1-a)}} \notag\\
                                                    &=  \sqrt{\frac{2}{\pi}}(1- \sqrt{1-t}-\sqrt{t}). 
\end{align}
Then we get easily \eqref{33} by adding \eqref{34} and \eqref{35}.  $\square$\\\\
\textbf{Remark:} The result can also be obtained by observing the duality $(1-A;V(B)_{1-t}-V(B)_1,t \in [0,1]) \stackrel{d}{=}(A;V(B)_t, t \in [0,1])$. \\

From Lemma \ref{pre1} and Lemma \ref{pre2} follows easily \eqref{23}.
\subsubsection{Variance for $V(B)$} 
We now turn to calculate $\mathbb{E}V(B)_t^2$ for $0 \leq t \leq 1$. We have
$$\mathbb{E}V(B)_t^2=\mathbb{E}(V(B)_t^2 1_{A>t})+ \mathbb{E}(V(B)_t^2 1_{A \leq t}).$$
\begin{lemma}
\label{pre}
\begin{equation}
\label{36}
\forall t \in [0,1],~\mathbb{E}(V(B)_t^2 1_{A>t})=3t-\frac{6t}{\pi} \arcsin\sqrt{t}-\frac{2}{\pi} \sqrt{t^3(1-t)} .
\end{equation}
\end{lemma}
\textbf{Proof:} According to the formula for $\mathbb{E}(B^{me})^2(t), t\in [0,1]$, we have:
\begin{align*}
\mathbb{E}(V(B)_t^2 1_{A>t})&=\int_t^1 \frac{a \mathbb{E}(B^{me})^2(\frac{t}{a})}{\pi \sqrt{a(1-a)}} da \\
                                              &\stackrel{\eqref{26}}{=}\frac{3t}{\pi} \int_t^1 \frac{1}{\sqrt{a(1-a)}}da -\frac{t^2}{\pi} \int_t^1 \frac{1}{\sqrt{a^3(1-a)}}da  \\
                                              &= \frac{3t}{\pi} \left[2 \arcsin \sqrt{a}\right]_t^1-\frac{t^2}{\pi} \left[-2 \sqrt{\frac{1-a}{a}}\right]_t^1=\eqref{36}.~~\square                                             
\end{align*} 
\begin{lemma}
\begin{equation}
\label{37}
\forall t \in [0,1],~\mathbb{E}(V(B)_t^2 1_{A \leq t})=\frac{4-2t}{\pi} (\arcsin \sqrt{t} -\sqrt{t(1-t)}).
\end{equation}
\end{lemma}
\textbf{Proof:} We have 
\begin{align}
\mathbb{E}(V(B)_t^2 1_{A \leq t}) &=\mathbb{E}((V(B)_t-V(B)_1)^2 1_{A \leq t})+\mathbb{E}(V(B)_1^2 1_{A \leq t}) \notag\\
                                                         &+2\mathbb{E}((V(B)_t-V(B)_1)V(B)_1 1_{A \leq t}). \notag
\end{align}
Denote $\beta_1, \beta_2$ and $\beta_3$ the three terms on the right hand side of the above equation. Note that $\beta_1$ can be easily computed from Lemma \ref{pre}  by change of variables:
\begin{equation}
\label{38}
\beta_1 = \frac{2(1-t)}{\pi}\left(3 \arcsin\sqrt{t}-\sqrt{t(1-t)}\right).  
\end{equation}
According to Denisov's decomposition, for $0 \leq a \leq t$, we have:
\begin{align}
\label{39}
\beta_2 &=\int_0^t \mathbb{E} (V(B)_1^2|A=a) \frac{da}{\pi \sqrt{a(1-a)}} \notag\\
       &=\int_0^t \left(2a+2(1-a)-2 \sqrt{\frac{\pi}{2} a} \sqrt{\frac{\pi}{2} (1-a)}\right) \frac{da}{\pi \sqrt{a(1-a)}} \notag\\
       &=\frac{4 \arcsin \sqrt{t}}{\pi}-t. 
 \end{align}      
 and 
$$ \beta_3 =\int_0^t (\gamma_1+\gamma_2)\frac{da}{\pi \sqrt{a(1-a)}},$$
 where
\begin{align}
\gamma_1&:=\mathbb{E}\left((V(B)_t-V(B)_1)V(B)_a|A=a\right) \notag\\
       &= \sqrt{\frac{\pi}{2}a} \times \sqrt{1-a} \sqrt{\frac{2}{\pi}} \left(\sqrt{\frac{1-t}{1-a}(1-\frac{1-t}{1-a})} +\arcsin \sqrt{\frac{1-t}{1-a}}\right) \notag\\
       &=\sqrt{a(1-a)} \left(\sqrt{\frac{1-t}{1-a}(1-\frac{1-t}{1-a})} +\arcsin \sqrt{\frac{1-t}{1-a}}\right). \notag
\end{align}
and
\begin{align}
 \gamma_2&:=\mathbb{E}\left((V(B)_t-V(B)_1)(V(B)_1-V(B)_a)|A=a\right) \notag\\
       &=-(1-a) \times 2 \sqrt{\frac{1-t}{1-a}}=-2 \sqrt{(1-t)(1-a)}. \notag
\end{align}
computed by Proposition \ref{compmean},
Therefore,
$$\int_0^t \gamma_1 \frac{da}{\pi \sqrt{a(1-a)}} = \frac{1}{\pi} \int_{1-t}^1  \left(\sqrt{\frac{1-t}{a}(1-\frac{1-t}{a})} +\arcsin \sqrt{\frac{1-t}{a}}\right) da.$$
Observe that 
\begin{align}
\label{40}
\int_{1-t}^1 \sqrt{\frac{1-t}{a}(1-\frac{1-t}{a})}da &=\left[2 \sqrt{(a-1+t)(1-t)}-2(1-t) \arcsin \sqrt{\frac{a-1+t}{a}}\right]_{1-t}^1 \notag\\
                                                                           &=2\sqrt{t(1-t)} - 2(1-t) \arcsin\sqrt{t}. 
\end{align}
and
\begin{align}
\label{41}
\int_{1-t}^1  \arcsin \sqrt{\frac{1-t}{a}}da &=\left[a \arcsin \sqrt{\frac{1-t}{a}}\right]_{1-t}^1-\int_{1-t}^1 a \left(\arcsin \sqrt{\frac{1-t}{a}}\right)^{'} da \notag\\
                                                                      &=\arcsin \sqrt{1-t} -\frac{\pi}{2}(1-t) +\frac{\sqrt{1-t}}{2} \int_{1-t}^1 \frac{da}{\sqrt{a-1+t} } \notag\\
                                                                      &=\arcsin \sqrt{1-t} -\frac{\pi}{2}(1-t) +\sqrt{t(1-t)}. 
\end{align}
We have by \eqref{40} and \eqref{41}:
\begin{equation}
\label{42}
\int_0^t \gamma_1 \frac{da}{\pi \sqrt{a(1-a)}}=\frac{t}{2}+\frac{2t-3}{\pi} \arcsin \sqrt{t} +\frac{3}{\pi} \sqrt{t(1-t)}.
\end{equation}
In addition,
\begin{equation}
\label{43}
\int_0^t \gamma_2 \frac{da}{\pi \sqrt{a(1-a)}}= -\frac{4}{\pi}\sqrt{t(1-t)}.
\end{equation}
Combining \eqref{38}, \eqref{39}, \eqref{42} and \eqref{43}, we obtain \eqref{37}.  $\square$\\

Finally, by adding \eqref{36} to \eqref{37}, we get \eqref{24}.
\section{Appendix: Computations for Brownian meander}
Let $(B^{me}_t, t \in [0,1])$ be standard Brownian meander. From Chung \cite{Chung}, we derive easily the density for meander along the paths as well as its joint distribution with terminal value:
\begin{equation}
\label{44}
\mathbb{P}(B^{me}_t \in dx)=t^{-\frac{3}{2}}x \exp\left(-\frac{x^2}{2t}\right) \erf\left(\frac{x}{\sqrt{2(1-t)}}\right)dx.
\end{equation}
\begin{equation}
\label{45}
\mathbb{P}(B^{me}_t \in dx, B^{me}_1 \in dy)=t^{-\frac{3}{2}}x \exp\left(-\frac{x^2}{2t}\right) \left( p_{1-t}(x,y)-
p_{1-t}(x,-y)\right) dydx.
\end{equation}
where $\erf$ is the error function for standard normal distribution: $\erf(x)=\frac{2}{\sqrt{\pi}} \int_{- \infty}^x \exp(-t^2)dt$ and $p$ is the transition kernel associated to Brownian motion: $p_{t}(x,y)=\frac{1}{\sqrt{2 \pi t}} \exp\left({-\frac{(x-y)^2}{2t}}\right)$.\\\\
\textbf{Proof of Proposition \ref{compmean}:}
\textbf{(a).} We compute the expectation of standard Brownian meander along the path, which relies on the following identity found in Gradshteyn and Ryzhik \cite{GR}:
\begin{equation}
\label{46}
\forall a>0,~\int_0^{\infty} x^2 \exp(-ax^2) \erf(x)dx =\frac{\sqrt{a}+(a+1)\arcsin\sqrt{\frac{1}{a+1}}}{2 \sqrt{\pi} a^{\frac{3}{2}}(a+1)} .
\end{equation}
By change of variables, we obtain:
\begin{align}
\mathbb{E}B^{me}_t&~~~~\stackrel{\eqref{44}}{=}\int_0^{\infty}   t^{-\frac{3}{2}}y^2 \exp\left(-\frac{y^2}{2t}\right) \erf\left(\frac{y}{\sqrt{2(1-t)}}\right) dy \notag\\
                                   &\stackrel{x=\frac{y}{\sqrt{2(1-t)}}}{=}\sqrt{8} \left(\frac{1-t}{t}\right)^{\frac{3}{2}} \int_0^{\infty} x^2 \exp\left(-\frac{1-t}{t}x^2\right) \erf(x) dx  \notag\\
                                   &~~~~\stackrel{\eqref{46}}{=}\sqrt{\frac{2}{\pi}} \left(\sqrt{t(1-t)}+ \arcsin \sqrt{t}\right). \notag
\end{align}
\textbf{(b).} We next calculate meander's second moment along the paths with the following identity also found in Gradshteyn and Ryzhik \cite{GR}: 
\begin{equation}
\label{47}
\forall a>0,~\int_0^{\infty} x^3 \exp(-ax^2) \erf(x)dx =\frac{2+3a}{4a^2(a+1)^{\frac{3}{2}}}.
\end{equation}
By change of variables, we get:
\begin{align}
\mathbb{E}(B^{me}_t)^2&~~~~\stackrel{\eqref{44}}{=}\int_0^{\infty}   t^{-\frac{3}{2}}y^3 \exp\left(-\frac{y^2}{2t}\right) \erf\left(\frac{y}{\sqrt{2(1-t)}}\right) dy \notag\\
                                           &\stackrel{x=\frac{y}{\sqrt{2(1-t)}}}{=}4t^{-\frac{3}{2}}(1-t)^2 \int_0^{\infty} x^3 \exp\left(-\frac{1-t}{t}x^2\right) \erf(x) dx  \notag\\
                                           &~~~~\stackrel{\eqref{47}}{=}  3t-t^2. \notag
\end{align}
\textbf{(c).} Finally we will compute $\mathbb{E}B^{me}_tB^{me}_1$ for $0 \leq t \leq 1$.
$$\mathbb{E}B^{me}_tB^{me}_1\stackrel{\eqref{45})}{=}\int_0^{\infty} \left(\int_0^{\infty} y(p_{1-t}(x,y)-
p_{1-t}(x,-y))dy\right) t^{-\frac{3}{2}}x^2 \exp\left(-\frac{x^2}{2t}\right)dx.$$
Remark that $\int_0^t y p_{t}(x,y) dy=\sqrt{\frac{t}{2 \pi}} e^{-\frac{x^2}{2t}}+x \erf(-\frac{x}{2 \sqrt{t}})$, we have:
$$\int_0^{\infty} y(p_{1-t}(x,y)-p_{1-t}(x,-y))dy=x \left(\erf(-\frac{x}{2 \sqrt{1-t}})+\erf(\frac{x}{2 \sqrt{1-t}})\right)=x.$$
Since it's well-known that for $a>0$, $\int_0^{\infty} x^3 e^{-ax^2} dx=\frac{1}{2a^2}$, we get:
$$\mathbb{E}B^{me}_tB^{me}_1=t^{-\frac{3}{2}} \int_0^{\infty} x^3 e^{-\frac{x^2}{2t}}dx=2 \sqrt{t}$$
\textbf{Remark:} The result in $(c)$, i.e. the identity \eqref{27} can be directly derived from Imhof's relation \cite{Imhof} (between Brownian meanders and three dimensional Bessel processes).
\\\\
\textbf{Acknowledgement:} The authors would like to express their gratitude to N.Forman for helpful discussion and suggestions throughout the preparation of this work. They would also like to thank P.Fitzsimmons for his remarks on the path decomposition result.
\addcontentsline{toc}{section}{References}
\bibliographystyle{plain}
\bibliography{Vervaat}
\end{document}